\documentclass[12pt]{article}

\usepackage{amsbsy}
\usepackage{amsmath}
\usepackage{amssymb}
\usepackage{epsfig}
\usepackage{multicol}
\usepackage{multirow}
\usepackage{color}
\usepackage{bm}
\usepackage{threeparttable}
\usepackage{float}
\renewcommand{\baselinestretch} {1.3}
\makeatletter 
\def\singlespace{\def\baselinestretch{1}\@normalsize}

\@addtoreset{equation}{section}

\newcommand{\qed}{$\Box$}

\newcommand{\bg}{\begin{eqnarray}}
\newcommand{\ed}{\end{eqnarray}}
\newcommand{\bgn}{\begin{eqnarray*}}
\newcommand{\edn}{\end{eqnarray*}}

\makeatletter
\def\singlespace{\def\baselinestretch{1}\@normalsize}

\date{\today}

\title{The consistency of estimator under fixed design regression model with NQD errors}

\author{
Jianhua Shi$^{1,2*}$, Xiaoping Chen$^{2,3}$, and Yong Zhou$^{2,4}$\\
 {\footnotesize $^{1}$School of Mathematics and Statistics,}
          {\footnotesize Minnan Normal University, Zhangzhou, China }\\
 {\footnotesize $^{2}$School of Statistics and Management,}
          {\footnotesize Shanghai University of Finance and Economics, Shanghai, China }\\
   {\footnotesize $^{3}$School of Mathematics and Computer Science,}
          {\footnotesize Fujian Normal University, Fuzhou, China }\\
   {\footnotesize $^{4}$Academy of Mathematics and Systems Science,
          Chinese Academy of Sciences, Beijing, China}\\
}

\begin{document}
\maketitle

\begin{singlespace}
\begin{footnotetext}
{\textbf{Acknowledgements}: The work was partially supported by National Natural Science Foundation of China (NSFC) (No.71271128,No.11301473), Science Fund for Creative Research Groups (11021161),Graduate Innovation Foundation of Shanghai University of Finance and Economics, China (CXJJ2012-423, CXJJ2013-451),and Natural Science Foundation of Fujian Province, China (2012J01028).

\textbf{*Corresponding author:}
Jianhua Shi, School of Statistics and Management, Shanghai University of Finance and Economics,Shanghai 200433, PR China, (\textit{E-mail: v0085@126.com}).
}
\end{footnotetext}
\end{singlespace}

\begin{abstract}
In this article, basing on NQD samples, we investigate the fixed design nonparametric regression model, i.e. $Y_{nk}  = g(x_{nk} ) + \varepsilon _{nk}$ for $1 \le k \le n$, where $\varepsilon _{nk}$ are pairwise NQD random errors, $x_{nk}$ are fixed design points, and $g( \cdot )$ is an unknown function. Nonparametric weighted estimator $g_n ( \cdot )$ of $g( \cdot )$ will be introduced and its consistency is studied. As special case, the consistency result for weighted kernel estimators of the model is established. This extends the earlier work on independent random and dependent random errors to NQD case.

\end{abstract}


\bigskip

\baselineskip=20pt
\section{Introduction}

 In regression analysis, it is common practice to investigate the functional relationship between the responses and design points. Nonparametric regression model provides a useful explanatory and diagnostic tool for this purpose. One may see Muller [1] and Hardle [2] for many examples about this and good introductions to the general subject area.

  To begin with, consider the fixed design nonparametric regression model in the paper
\[
Y_{nk}  = g(x_{nk} ) + \varepsilon _{nk} ,_{} 1 \le k \le n.
\]
Here $x_{nk} ,1 \le k \le n$ , are known fixed design points, and $\varepsilon _{nk}$  are random errors , $g( \cdot )$ is an unknown regression function. As an estimate of $g( \cdot )$, we consider the following general linear smoother.
\[
g_n (x) = \sum\limits_{k = 1}^n {\omega _{nk} (x)Y_{nk} },
\]
where the weight functions $\omega _{nk} (x)$ depend on $x,x_{n1} , \cdots ,x_{nn}$.

It is well known that Georgiev [3] first proposed the estimator above, and the estimator subsequently have been studied by many authors. A brief review of the theoretic development in recent years is worth mentioning. Results on $\varepsilon _{nk}$ being assumed to be independent, consistency and asymptotic normality have been investigated by Georgiev [4] and M\"uller [5] among others. Results for the case when $\varepsilon _{nk}$ are dependent have also been studied by various authors in recent years. Roussas et al. [6] established asymptotic normality of $g_n (x)$ assuming that the errors are from a strictly stationary stochastic process under the strong mixing condition. Tran et al. [7] discussed again asymptotic normality of $g_n (x)$ assuming that the errors form a weakly stationary linear process with a martingale difference sequence. Hu et al. [8] gave the mean consistency, complete consistency, and asymptotic normality of regression models based on linear process errors. Under negatively associated sequences, Liang and Jing [9] presented some asymptotic properties for estimates of nonparametric regression models, Yang et al. [10] generalized part results of Liang and Jing [9] for negatively associated sequences to the case of negatively orthant dependent sequences, and so on.

In this paper, we shall investigate the above nonparametric regression problem under pairwise NQD errors, which means more general case for sampling.

\textbf{Definition 1.1.} [11] The pair $(X,Y)$
 of random variables X and Y is said to be NQD (negatively quadrant dependent), if
\[
               P(X \le x,Y \le y) \le P(X \le x)P(Y \le y),\forall x,y \in R.             \eqno (1.1)
\]
A sequence of random variables $\{ X_n ,n \ge 1\} $ is pairwise NQD random variables(for short, NQD), if $(X_i ,Y_j )$ is NQD for every   $i \ne j,i,j = 1,2, \cdots $ .

It can be deduced from Definition 1.1 that
\[
         P(X > x,Y > y) \le P(X > x)P(Y > y),\forall x,y \in R.              \eqno (1.2)
\]
Moreover, it follows that (1.2) also implies (1.1), and hence, (1.1) and (1.2) are actually equivalent.

The definition was introduced by Lehmann [11], which contains independent random variable, NA(negatively associated) random variable and NOD(negatively orthant dependent) random variable et al. as special cases. For the reason of the wide applications of NQD random variables in reliability theory and applications, the notions of NQD random variables have received many  concern recently. Some properties about NQD random variables can be found in Lehmann [4], and there are many other meaningful literature (e.g. Matula [12], Huang et al. [13], Sung [14], Shi [15], Wang et al. [17], Li et al. [18]).

However, the pairwise NQD structure is more comprehensive than the NA (negative associated) structure and the NOD (negatively orthant dependent) structure. Concerning to the study for the theory of pairwise NQD random variables, due to lack of some key technique tool, such as Bernstain type inequality and exponential inequality etc. still unestablished for NQD sequences, investigating related result is restraint, especially the estimators of parametric and nonparametric components in regressions model under NQD error's structure. Hence, extending the asymptotic properties of independent and other dependent random variables to the case of NQD variables is highly desirable and of considerably significance in the theory and application.

In this article, basing on several related lemmas, we investigate the fixed design nonparametric regression model with NQD errors. Nonparametric estimator $g_n ( \cdot )$ of $g( \cdot )$  will be introduced and its usual consistency properties of $g_n ( \cdot )$ including mean convergence, uniform mean convergence, convergence in probability, et al. are studied under suitable regularity conditions.

The organization of this paper is as follows. In section 2, we shall present several lemmas for proof of main results, and give the basic assumptions for the nonparametric estimator. We give the further assumption and the main results in section 3. The proofs of the results will be deferred to Section 4.

\section{Some lemmas and Basic assumptions} 
\subsection{ Some lemmas}
We shall begin with a few preliminary lemmas useful in the proofs of our main results. Firstly, a fact about the NQD properties is cited from [11].\\
\textbf{Lemma 2.1.}
\emph{[11] Let the pair $(X,Y)$ of random variables $X$ and $Y$ be NQD, then}

\emph{(1) $E(XY) \le EX \cdot EY$ ;}

\emph{(2)  $P(X > x,Y > y) \le P(X > x)P(Y > y),$ for any $x,y \in R$;}

\emph{ (3) If $f,g$ are both non-decreasing (or non-increasing) functions, then $f(X)$ and $g(X)$ are NQD.}\\
\textbf{Lemma 2.2.}
\emph{[16] Let $\{ X_n ,n \ge 1\}$ be a sequence of pairwise NQD random variables such that $EX_n  = 0,EX_n^2  < \infty$ for all $n \ge 1$, denote  $T_j (k) \buildrel \Delta \over = \sum\nolimits_{i = j + 1}^{j + k} {X_i } ,j \ge 0,k \ge 1$ ,\ then
 \[
E(T_j (k))^2  \le \sum\limits_{i = j + 1}^{j + k} {EX_i^2 } ,E\mathop {\max }\limits_{1 \le k \le n} (T_j (k))^2  \le \frac{{4\log ^2 n}}{{\log ^2 2}}\sum\limits_{i = j + 1}^{j + k} {EX_i^2 } .
\]}
\ \ In the rest below, we assume $0 = x_{n({\rm{0)}}}  \le x_{n({\rm{1)}}}  \le x_{n({\rm{2)}}}  \le  \cdots  \le x_{n(n)}  = 1$ and let $
\delta _n  = \mathop {\max }\limits_{1 \le k \le n} (x_{n(k)}  - x_{n(k - 1)} )$ . Furthermore, assume that

$(A_1 )$  $K( \cdot )$ is bounded and satisfies Lipschitz condition of order $\alpha (\alpha  > 0)$ on $R^1$ , and $\int_{ - \infty }^\infty  {\left| {K(u)} \right|du}  < \infty$ ;

$(A_2 )$  $h_n  \to 0$  and ${{\delta _n^\alpha  } \mathord{\left/
 {\vphantom {{\delta _n^\alpha  } {h_n^{1 + \alpha } }}} \right.
 \kern-\nulldelimiterspace} {h_n^{1 + \alpha } }} \to 0$ as $n \to \infty$;

$(A_3 )$ $\int_{ - \infty }^\infty  {K(u)du}  = 1$ ;\\
\textbf{Lemma 2.3.}
\emph{If Conditions $(A_1 ),(A_2 )$ hold, then
 \[
\mathop {\lim }\limits_{n \to \infty } \sum\limits_{k = 1}^n {\frac{{x_{n(k)}  - x_{n(k - 1)} }}{{h_n }}\left| {K(\frac{{x - x_{n(k)} }}{{h_n }})} \right|}  = \int_{ - \infty }^{ + \infty } {\left| {K(u)} \right|du} ,_{} x \in (0,1), \eqno (2.1)
\]
and for a fixed point $\tau  \in (0,{1 \mathord{\left/
 {\vphantom {1 2}} \right.
 \kern-\nulldelimiterspace} 2})$,
 \[
\mathop {\lim }\limits_{n \to \infty } \mathop {\sup }\limits_{x \in [\tau ,1 - \tau ]} \sum\limits_{k = 1}^n {\frac{{x_{n(k)}  - x_{n(k - 1)} }}{{h_n }}\left| {K(\frac{{x - x_{n(k)} }}{{h_n }})} \right|}  = \int_{ - \infty }^{ + \infty } {\left| {K(u)} \right|du}.  \eqno (2.2)
\]}
\emph{Proof of Lemma 2.3} \ Denote $H(x) = I(0 \le x \le 1)$, where $I( \cdot )$ is the usual indicator function, and
\begin{eqnarray*}
 &&\sum\limits_{k = 1}^n {\frac{{x_{n(k)}  - x_{n(k - 1)} }}{{h_n }}\left| {K(\frac{{x - x_{n(k)} }}{{h_n }})} \right|}  - \int_{ - \infty }^{ + \infty } {\left| {K(u)} \right|du}  \\
 & =& \{ \sum\limits_{k = 1}^n {\frac{{x_{n(k)}  - x_{n(k - 1)} }}{{h_n }}\left| {K(\frac{{x - x_{n(k)} }}{{h_n }})H(x_{n(k)} )} \right|}  - h_n^{ - 1} \int_0^1 {\left| {K(\frac{{x - u}}{{h_n }})} \right|H(u)du\} }  \\
 && + \{ h_n^{ - 1} \int_{ - \infty }^{ + \infty } {\left| {K(\frac{{x - u}}{{h_n }})} \right|H(u)du}  - \int_{ - \infty }^{ + \infty } {\left| {K(u)} \right|du\} }  \buildrel \Delta \over = T_{n1} (x) + T_{n2} (x)
\end{eqnarray*}

We look at each term separately. Note that there is $\theta _{n(k)}  \in (0,1),^{} (k = 1,2, \cdots ,n)$ by Mean- value Theorem for integrals such that
\begin{eqnarray*}
 \left| {T_{n1} (x)} \right| &=& \left| {h_n^{ - 1} \sum\limits_{k = 1}^n {\tilde \delta _{n(k)} \{ \left| {K(\frac{{x - x_{n(k)} }}{{h_n }})} \right|H(x_{n(k)} )} } \right. \\
& & \left. { - \left| {K(\frac{{x - x_{n(k)}  + \theta _{n(k)} \tilde \delta _{n(k)} }}{{h_n }})} \right|H(x_{n(k)}  - \theta _{n(k)} \tilde \delta _{n(k)} )\} } \right| \\
  &\le& h_n^{ - 1} \sum\limits_{k = 1}^n {\tilde \delta _{n(k)} \{ \left| {K(\frac{{x - x_{n(k)} }}{{h_n }}) - K(\frac{{x - x_{n(k)}  + \theta _{n(k)} \tilde \delta _{n,k} }}{{h_n }})} \right|\left| {H(x_{n(k)} )} \right|}  \\
 & & + \left| {K(\frac{{x - x_{n(k)}  + \theta _{n(k)} \tilde \delta _{n(k)} }}{{h_n }})} \right|\left| {H(x_{n(k)} ) - H(x_{n(k)}  - \theta _{n(k)} \tilde \delta _{n(k)} )} \right|\}  \\
 & \le& Mh_n^{ - 1} \sum\limits_{k = 1}^n {\tilde \delta _{n(k)} ({{\theta _{n(k)} \tilde \delta _{n(k)} } \mathord{\left/
 {\vphantom {{\theta _{n(k)} \tilde \delta _{n(k)} } {h_n }}} \right.
 \kern-\nulldelimiterspace} {h_n }})^\alpha  }  \le M{{({{\delta _n } \mathord{\left/
 {\vphantom {{\delta _n } {h_n }}} \right.
 \kern-\nulldelimiterspace} {h_n }})^\alpha  } \mathord{\left/
 {\vphantom {{({{\delta _n } \mathord{\left/
 {\vphantom {{\delta _n } {h_n }}} \right.
 \kern-\nulldelimiterspace} {h_n }})^\alpha  } {h_n }}} \right.
 \kern-\nulldelimiterspace} {h_n }},
\end{eqnarray*}
where $\tilde \delta _{n(k)}  = x_{n(k)}  - x_{n(k - 1)}$ , and we use Condition $(A_1 )$ to the second inequality .

Then, according to Condition $(A_2 )$, we conclude
\[
\mathop {\lim }\limits_{n \to \infty } \left| {T_{n1} (x)} \right| = 0, \
and  \ \
\mathop {\lim }\limits_{n \to \infty } \mathop {\sup }\limits_{x \in [\tau ,1 - \tau ]} \left| {T_{n1} (x)} \right| = 0.      \eqno (2.3)
\]

As for $T_{n2} (x)$, when $x \in (0,1)$, we have
\[
\left| {T_{n2} (x)} \right| \le \int_{ - \infty }^{ + \infty } {\left| {H(x - h_n u) - H(x)} \right|\left| {K(u)} \right|du} .
\]

Note that by the definition of $H( \cdot )$, $\mathop {\lim }\limits_{n \to \infty } H(x - h_n u) = H(x)$ for all $x \in (0,1)$ and $u \in R^1$. Under the integrability of $\left| {K(u)} \right|$, $\left| {T_{n2} (x)} \right| \to 0,n \to \infty$ \ by Dominated Convergence Theorem, which together with (2.3) implies (2.1).

Again because of $\int_{ - \infty }^\infty  {\left| {K(u)} \right|du < \infty }$ and $h_n  \to 0$, if $n$ large sufficiently, one can choose a sufficient small positive number $\tau _0$, such that when $\left| {h_n u} \right| < \tau _0  < \tau$, there is
 \[
\int_{ - \infty }^\infty  {\left| {K(u)} \right|I(\left| u \right| \ge {{\tau _0 } \mathord{\left/
 {\vphantom {{\tau _0 } {h_n }}} \right.
 \kern-\nulldelimiterspace} {h_n }})du}  < \frac{\varepsilon }{2}.
\]

As a result, for $x \in [\tau ,1 - \tau ]$, uniformly
\begin{eqnarray*}
\left| {T_{n2} (x)} \right| &\le& \int_{ - \infty }^{ + \infty } {\left| {H(x - h_n u) - H(x)} \right|\left| {K(u)} \right|[I(\left| {h_n u} \right| < \tau _0 ) + I(\left| {h_n u} \right| \ge \tau _0 )]du}\\
& \le& 2\int_{ - \infty }^{ + \infty } {\left| {K(u)} \right|I(\left| {h_n u} \right| \ge \tau _0 )du}  < \varepsilon
\end{eqnarray*}

Therefore,
 \[
\mathop {\sup }\limits_{x \in [\tau ,1 - \tau ]} \left| {T_{n2} (x)} \right| \to 0{\rm{,}}_{}\ as \ n \to \infty .
\]

Combining (2.3), then (2.2) holds, as we wanted to show. This completes the proof. \qed \\
\textbf{Lemma 2.4.}
\emph{If Conditions $(A_1 ),(A_2 ),(A_3 )$ hold, then
\[
\mathop {\lim }\limits_{n \to \infty } \sum\limits_{k = 1}^n {\frac{{x_{n(k)}  - x_{n(k - 1)} }}{{h_n }}K(\frac{{x - x_{n(k)} }}{{h_n }})}  = 1,_{} x \in (0,1), \eqno (2.4)
\]
and for a fixed point$ \ \tau  \in (0,{1 \mathord{\left/ {\vphantom {1 2}} \right. \kern-\nulldelimiterspace} 2}),$
\[
\mathop {\lim }\limits_{n \to \infty } \mathop {\sup }\limits_{x \in [\tau ,1 - \tau ]} \left| {\sum\limits_{k = 1}^n {\frac{{x_{n(k)}  - x_{n(k - 1)} }}{{h_n }}K(\frac{{x - x_{n(k)} }}{{h_n }})}  - 1} \right| = 0. \eqno (2.5)
\]}
\emph{Proof of Lemma 2.4} The proof is similar to those of Lemma 2.3 with $\left| {K(u)} \right|$ replaced by ${K(u)}$ and using Condition $(A_3 )$, so is omitted here. \ \ \ \ \ \ \ \ \ \ \ \ \ \ \ \ \ \ \ \ \ \ \ \ \ \ \ \ \ \ \ \ \ \ \ \ \ \ \ \ \ \ \ \ \ \ \ \ \ \ \qed

\subsection{Basic assumptions}

Unless otherwise specified, we assume throughout the paper that the random sample $(x_{nk} ,Y_{nk} )$ for $1 \le k \le n$ come from the regression model
 \[
Y_{nk}  = g(x_{nk} ) + \varepsilon _{nk} ,1 \le k \le n,\eqno (2.6)
\]
where $\{ \varepsilon _{nk} ,1 \le k \le n\}$ from a sequence of zero mean random errors with the same distribution as $\{ \varepsilon _k ,1 \le k \le n\}$ for each $n$, $\{ x_{nk} ,1 \le k \le n\}$ are known fixed design points from a compact set $A$ in $R^d$ ($d$  is a positive integer), and $g( \cdot )$ is an unknown real valued regression function and assumed to be bounded on the compact set $A$.

The present paper investigates the general linear smoother as an estimate of $g( \cdot )$ in the following, defined by formula
 \[
g_n (x) = \sum\limits_{k = 1}^n {\omega _{nk} (x)Y_{nk} } ,_{} 1 \le k \le n,
\]
where the array of weight functions $\omega _{nk} (x),1 \le k \le n$ depends on the fixed design points $x,x_{n1} , \cdots ,x_{nn}$ and on the number of observations $n$, which $\omega _{nk} (x) = 0$ for $k > n$.

In the following section, we denote all continuity points of the function $g( \cdot )$ on set $A$ as $C(g)$. Let the symbol $\left\| x \right\|$ be the Eucledean norm of $x$ , $M$ a generic positive constant in the sequel, which could take different values at different places.

\section{Main results}

We shall establish two different models of convergence for the nonparametric regression estimate $g_n (x)$ at a fixed point $x$. First, we give some assumptions on weight function $\omega _{nk} (x)$ in the following. The similar assumptions on weighted functions can be found in Georgiev et al. [4], Hu et al. [8], Liang et al. [9]and Yang et al. [10], etc.
\begin{description}
\item[\ \ \ \ ($B_1$ )] $\sum\limits_{k = 1}^n {\omega _{nk} (x)}  \to 1,\ as_{} \ n \to \infty ;$
 \item[\ \ \ \ ($B_2$ )] $\sum\limits_{k = 1}^n {\left| {\omega _{nk} (x)} \right|}  \le M,\ \forall n;$
   \item[\ \ \ \ ($B_3$ )] $\sum\limits_{k = 1}^n {\omega _{nk}^2 (x)}  \to 0,\ as_{}\  n \to \infty ;$
   \item[\ \ \ \ ($B_4$ )] $\sum\limits_{k = 1}^n {\left| {\omega _{nk} (x)} \right|} I(\left\| {x_{nk}  - x} \right\| > a) \to 0,\ as_{}\  n \to \infty ,for_{}\  a > 0.$
\end{description}

The weights $\omega _{nk} (x),1 \le k \le n$ in the assumptions is relatively extensive in practice, which can be easily satisfied by the commonly adopted weights used, such as the well-known nearest neighbor weights.\\
\textbf{Example 3.1.} Let $g( \cdot )$ be continuous on interval $A \buildrel \Delta \over = [0,1]$. Without loss of generality, put $x_{nk}  = {k \mathord{\left/ {\vphantom {k n}} \right. \kern-\nulldelimiterspace} n},1 \le k \le n$. When $\left| {x_{ni}  - x} \right| = \left| {x_{nj}  - x} \right|$, assume that $\left| {x_{ni}  - x} \right|$ is ahead of $\left| {x_{nj}  - x} \right|$ for $x_{ni}  < x_{nj}$, then a permutation for $\left| {x_{n1}  - x} \right|,\left| {x_{n2}  - x} \right|, \cdots ,\left| {x_{nn}  - x} \right|$ can be given as follows
 \[
\left| {x_{R_1 (x)}^n  - x} \right| \le \left| {x_{R_2 (x)}^n  - x} \right| \le  \cdots  \le \left| {x_{R_n (x)}^n  - x} \right|,x \in A.
\]

Let $k_n  = o(n)$, if define the nearest neighbor weight as
\[
\omega _{nk} (x) = \left\{ {\begin{array}{*{20}c}
   {{1 \mathord{\left/
 {\vphantom {1 {k_n ,}}} \right.
 \kern-\nulldelimiterspace} {k_n ,}}} & {\left| {x_{nk}  - x} \right| \le \left| {x_{R_{k_n } (x)}^n  - x} \right|,}  \\
   {0,} & {otherwise.}  \\
\end{array}} \right.
\]
Then, one can easily verify by the choice of $x_{ni}$ and the definition of $R_i (x)$ that Conditions $(B_1 ) \sim (B_4 )$ are satisfied.

We now state our first result for the mean convergence of $g_n (x)$, which, on the opinion of statistics, is asymptotically unbiased of $g (x)$ in the proof of Theorem 3.1.\\
\textbf{Theorem 3.1.}\emph{(Mean convergence) Assume that Conditions $(B_1 ) \sim (B_4 )$ hold. Let \ $\{ \varepsilon _n ,n \ge 1\}$ be a mean zero pairwise NQD sequences with $\mathop {\sup }\limits_{n \ge 1} E\varepsilon _n^2  < \infty$, if \ $0 < p \le 2$, then
\[
\mathop {\lim }\limits_{n \to \infty } E(g_n (x) - g(x))^p  = 0,  \eqno (3.1)
\]
for $\forall x \in C(g)$.}

Another similar form of mean convergence,  by using the inequality $(\sum\nolimits_{k = 1}^n {\left| {a_i } \right|^\beta  } )^{{1 \mathord{\left/ {\vphantom {1 \beta }} \right. \kern-\nulldelimiterspace} \beta }}  \le (\sum\nolimits_{k = 1}^n {\left| {a_i } \right|^\alpha  } )^{{1 \mathord{\left/ {\vphantom {1 \alpha }} \right. \kern-\nulldelimiterspace} \alpha }} ,1 \le \alpha  \le \beta ,$ for any real number sequence $\{ a_i ,1 \le i \le n\}$, is that.\\
\textbf{Theorem ${\rm{3}}{\rm{.1.}}^\prime$}\emph{(Mean convergence) Assume that Conditions $(B_1 ),(B_2 ) \ and \  (B_4 )$ hold. Let \ $\{ \varepsilon _n ,n \ge 1\}$ be a mean zero pairwise NQD sequences with $\mathop {\sup }\limits_{n \ge 1} E\varepsilon _n^p  < \infty$ for some $1 < p \le 2$, if $\sum\limits_{k = 1}^n {\left| {\omega _{nk}^s (x)} \right|}  \to 0,as\ n \to \infty ,$ with $1 < s \le p$ . then (3.1) holds
for $\forall x \in C(g)$.}

For any fixed point $x$ on a compact set $A$ in $R^d ,(d \ge 1)$, in order to obtaining uniform convergence for the estimator of $g(x)$, some uniform version of assumptions on $\omega _{nk} (x)$ are necessarily replaced by that as follows.

\begin{description}
    \item[\ \ \ \ ($B'_1$ )] $\mathop {\sup }\limits_{x \in A} \left| {\sum\limits_{k = 1}^n {\omega _{nk} (x)}  - 1} \right| \to 0,\ as_{}\ n \to \infty ;$
 \item[\ \ \ \ ($B'_2$ )] $\mathop {\sup }\limits_{x \in A} \left| {\sum\limits_{k = 1}^n {\omega _{nk} (x)} } \right| \le M,\ \forall n;$
   \item[\ \ \ \ ($B'_3$ )] $\mathop {\sup }\limits_{x \in A} \sum\limits_{k = 1}^n {\omega _{nk}^2 (x)}  \to 0,\ as_{} \ n \to \infty ;$
   \item[\ \ \ \ ($B'_4$ )] $\mathop {\sup }\limits_{x \in A} \sum\limits_{k = 1}^n {\left| {\omega _{nk} (x)} \right|} I(\left\| {x_{nk}  - x} \right\| > a) \to 0,as_{} \ n \to \infty ,for_{} \ a > 0.$
\end{description}

Then we are in the position to give the following result.\\
\textbf{Theorem 3.2.}\emph{(uniform mean convergence) Assume that Conditions $(B'_1 ) \sim (B'_4 )$ hold. Let $g( \cdot )$ be continuous on the compact set $A$ , $\{ \varepsilon _n ,n \ge 1\}$ a mean zero NQD sequence. If $\mathop {\sup }\limits_{n \ge 1} E\varepsilon _n^2  < \infty$, then
\[
\mathop {\lim }\limits_{n \to \infty } \mathop {\sup }\limits_{x \in A} E(g_n (x) - g(x))^p  = 0,\eqno (3.2)
\]
for $0 < p \le 2$.}\\
\textbf{Remark 3.1.} Since NA sequence and NOD sequence are NQD sequence, we generalize some results of Liang et al. [9]and Yang et al. [10] to the case of NQD errors, respectively. And as a consequence, one may get consistency property for the weighted kernel estimators in the model (2.6).\\
\textbf{Corollary 3.1.} \emph{Assume that Conditions $(A_1 ),(A_2 ),(A_3 )$ hold, and
 \[
\ \ (A_4 )\ \ \sum\limits_{k = 1}^n {\frac{{x_{n(k)}  - x_{n(k - 1)} }}{{h_n }}\left| {K(\frac{{x - x_{n(k)} }}{{h_n }})} \right|I(\left| { x_{n(k)} -x } \right| > a)}  \to 0,_{}^{} as_{}\  n \to \infty ,_{}^{} for_{} \ a > 0.
\]
Let $\{ \varepsilon _n ,n \ge 1\}$ be a mean zero pairwise NQD sequences with $\mathop {\sup }\limits_{n \ge 1} E\varepsilon _n^2  < \infty
$, $g( \cdot )$ a continuous function on interval $(0,1)$, if $0 < p \le 2$, then
 \[
\mathop {\lim }\limits_{n \to \infty } E(\sum\limits_{k = 1}^n {\frac{{x_{n(k)}  - x_{n(k - 1)} }}{{h_n }}K(\frac{{x - x_{n(k)} }}{{h_n }})Y_{nk} }  - g(x))^p  = 0, \eqno (3.3)
\]
for $\forall x \in (0,1)$.}

\emph{Furthermore, if $(A_4 )$ is replaced by}

\emph{$(A_4 )'\mathop {\sup }\limits_{x \in [\tau ,1 - \tau ]} \sum\limits_{k = 1}^n {\frac{{x_{n(k)}  - x_{n(k - 1)} }}{{h_n }}\left| {K(\frac{{x - x_{n(k)} }}{{h_n }})} \right|I(\left| { x_{n(k)} -x } \right| > a)}  \to 0,$
for $a > 0,$ \\
then
 \[
\mathop {\lim }\limits_{n \to \infty } \mathop {\sup }\limits_{x \in [\tau ,1 - \tau ]} E(\sum\limits_{k = 1}^n {\frac{{x_{n(k)}  - x_{n(k - 1)} }}{{h_n }}K(\frac{{x - x_{n(k)} }}{{h_n }})Y_{nk} }  - g(x))^p  = 0, \eqno (3.4)
\]
for $0 < p \le 2$.}

Next, we shall give the weak consistency for the estimator of $g(x)$ under existence of absolute mean for variable.\\
\textbf{Theorem 3.3.} \emph{(convergence in probability) Assume that Conditions $(B_1 ),(B_{\rm{2}} ),(B_{\rm{4}} )$ hold. Let $\{ \varepsilon _n ,n \ge 1\}$ be a mean zero pairwise NQD sequences and uniformly bounded by a random variable $X$ in the sense that $\mathop {\sup }\limits_{n \ge 1} P(\left| {\varepsilon _n } \right| \ge x) \le P(\left| X \right| \ge x)$ for all $x > 0$. If $E\left| X \right| < \infty$ and $\mathop {\sup }\limits_k \left| {\omega _{nk} (x)} \right| = o(1)$ , then
\[
g_n (x) \to g(x)_{}^{}, \ in\ probability \ as \ n \to \infty ,    \eqno (3.5)
\]
for $\forall x \in C(g)$.}\\
\textbf{Corollary 3.2.} Assume that Conditions $(A_1 ) \sim (A_4 )$ hold, and $\{ \varepsilon _n ,n \ge 1\}$ be a mean zero pairwise NQD sequences and uniformly bounded by a random variable $X$ in the sense that $\mathop {\sup }\limits_{n \ge 1} P(\left| {\varepsilon _n } \right| \ge x) \le P(\left| X \right| \ge x)$ for all $x > 0$, $g( \cdot )$ a continuous function on interval $(0,1)$. If $E\left| X \right| < \infty$ and \ $\delta _n h_n^{ - 1}  = o(1)$, then
\[
\sum\limits_{k = 1}^n {\frac{{x_{n(k)}  - x_{n(k - 1)} }}{{h_n }}K(\frac{{x - x_{n(k)} }}{{h_n }})Y_{nk} }  \to g(x)_{}^{} , \ in\ probability \ as \ n \to \infty ,
\]
for $\forall x \in (0,1)$.

\section{Proofs for main results}

\emph{Proof of Theorem 3.1} We write firstly the triangle inequality that
\begin{eqnarray}\label{eq4.1}
E\left| {g_n (x) - g(x)} \right|^p  \le ME\{ \left| {g_n (x) - Eg_n (x)} \right|^p  + \left| {Eg_n (x) - g(x)} \right|^p \} .
\end{eqnarray}

By Jensen's inequality, Lemma 2.2, Condition $(B_3 )$ and $\mathop {\sup }\limits_{n \ge 1} E\varepsilon _n^2  < \infty$, and note that $\{ \varepsilon _{nk} ,1 \le k \le n\}$ has the same distribution as $\{ \varepsilon _k ,1 \le k \le n\}$, we have for $0 < p \le 2$
\begin{eqnarray*}
\ \ \ \ \ \ \ \ \ \ \ \ \ E\left| {g_n (x) - Eg_n (x)} \right|^p  &\le& (E(\sum\limits_{k = 1}^n {\omega _{nk}^{} (x)\varepsilon _k^{} } )^2 )^{{p \mathord{\left/
 {\vphantom {p 2}} \right.
 \kern-\nulldelimiterspace} 2}}\\
 &\le& (\sum\limits_{k = 1}^n {\omega _{nk}^2 (x)E\varepsilon _k^2 } )^{{p \mathord{\left/
 {\vphantom {p 2}} \right.
 \kern-\nulldelimiterspace} 2}}  \to 0,_{}^{} as \ n \to \infty , \ \ \ \ \ \ \ \ \ \ \   \ \ \ \ \ \ \ \ \ \ (4.2)
\end{eqnarray*}
$x \in C(g)$,\ since $\{ \omega _{nk}^{} (x)\varepsilon _k^{} ,1 \le k \le n\}$ are also NQD according to Lemma 2.1.

Meanwhile, for the bias $Eg_n (x) - g(x)$, choose a number $a>0$, we can get the following upper bound:
\begin{eqnarray*}
\left| {Eg_n (x) - g(x)} \right| &\le& \sum\limits_{k = 1}^n {\left| {\omega _{nk} (x)} \right| \cdot \left| {g(x_{nk} ) - g(x))} \right|[I(\left\| {x_{nk}  - x} \right\| \le a)}\\
&& + I(\left\| {x_{nk}  - x} \right\| > a)] + \left| {\sum\limits_{k = 1}^n {\omega _{nk} (x)}  - 1} \right| \cdot \left| {g(x)} \right|,x \in C(g).
\end{eqnarray*}

Because of $x \in C(g)$, for any $\varepsilon  > 0$, there exists a $\delta  > 0$ such that $\left| {g(x_{nk} ) - g(x)} \right| < \varepsilon$ whenever $\left\| {x_{nk}  - x} \right\| < \delta$. Thus, by setting $0 < a < \delta$, Conditions $(B_1 ),(B_2 ),(B_4 )$ together with the arbitrary of $\varepsilon  > 0$ imply that the estimate $g_n ( \cdot )$ is asymptotically unbiased for $g ( \cdot )$, and then
 \[
\left| {Eg_n (x) - g(x)} \right|^p  \to 0,_{}^{} n \to \infty ,_{}^{} x \in C(g).\eqno (4.3)
\]

Therefore, we can deduce from (4.1),(4.2),(4.3) that (3.1) follows, this ends the proof. \qed \\
\emph{Proof of Theorem 3.2} Note that in a compact set $A$, $g( \cdot )$ is uniformly continuous if it is continuous. Consequently, similar proof as Theorem 3.1, we can get that
\begin{eqnarray*}
&&\mathop {\sup }\limits_{x \in A} E(g_n (x) - g(x))^p\\
&\le& M\{ \mathop {\sup }\limits_{x \in A} E\left| {g_n (x) - Eg_n (x)} \right|^p  + \mathop {\sup }\limits_{x \in A} \left| {Eg_n (x) - g(x)} \right|^p \} ,
\end{eqnarray*}
tends to zero if $n \to \infty$, which means the desired result (3.2).\ \ \ \ \ \ \ \ \ \ \ \ \ \ \ \ \ \ \ \ \ \ \ \ \ \ \ \ \ \ \ \ \ \ \ \ \ \ \ \ \ \qed \\
\emph{Proof of Corollary 3.1} Note that under Condition $(A_2 )$,
\[
\sum\limits_{k = 1}^n {(\frac{{x_{n(k)}  - x_{n(k - 1)} }}{{h_n }}K(\frac{{x - x_{n(k)} }}{{h_n }}))^2 }  \le M\frac{{\delta _n }}{{h_n }}\sum\limits_{k = 1}^n {\left| {\frac{{x_{n(k)}  - x_{n(k - 1)} }}{{h_n }}K(\frac{{x - x_{n(k)} }}{{h_n }})} \right|}  \to 0,
\]
as $n \to \infty$.

And there are
\[
\sum\limits_{k = 1}^n {\frac{{x_{n(k)}  - x_{n(k - 1)} }}{{h_n }}\left| {K(\frac{{x - x_{n(k)} }}{{h_n }})} \right|}  \le M,\forall n,
\]
\[
\sum\limits_{k = 1}^n {\frac{{x_{n(k)}  - x_{n(k - 1)} }}{{h_n }}K(\frac{{x - x_{n(k)} }}{{h_n }})}  \to 1,x \in (0,1),
\]
by Lemma2.3 and Lemma2.4, respectively. Therefore, according to Theorem 3.1, (3.3) follows.

As to (3.4), similar as above,  one may verify, on the interval $[\tau ,1 - \tau ]$,  the condition of Theorem 3.2 by the second result of Lemma2.3 and Lemma2.4, i.e. (2.2) and (2.5), respectively. This ends the proof.\ \ \ \ \ \ \ \ \ \ \ \ \ \ \ \ \ \ \ \ \ \ \ \ \ \ \ \ \ \ \ \ \ \ \ \ \ \ \ \ \ \ \ \ \ \ \ \ \ \ \ \ \ \ \ \ \ \ \ \ \ \ \ \ \ \ \ \ \ \ \ \ \ \ \ \ \ \ \ \ \ \ \ \ \ \ \ \ \ \ \ \ \ \ \ \ \ \ \ \ \ \ \ \ \ \ \ \ \ \ \ \ \ \ \ \ \ \ \ \ \ \ \ \ \ \ \qed \\
\emph{Proof of Theorem 3.3} Since $x \in C(g)$, the same reason as before, for any $\varepsilon  > 0$, there is a number $\delta  > 0$, when $a \in (0,\delta )$, one may get that $\left| {Eg_n (x) - g(x)} \right|$ tends to zero by arbitrary of $\varepsilon  > 0$
and Conditions $(B_1 ),(B_2 ),(B_4 )$. For proving (3.5), note that
\[
\left| {g_n (x) - g(x)} \right| \le \left| {g_n (x) - Eg_n (x)} \right| + \left| {Eg_n (x) - g(x)} \right|  \eqno (4.4)
\]

We now prove that the random part of r.h.s. in (4.4) tends to zero in probability as $n \to \infty$ . Observe that
\[
\left| {g_n (x) - Eg_n (x)} \right| = \left| {\sum\limits_{k = 1}^n {\omega _{nk} (x)\varepsilon _{nk} } } \right|
\]

Next introduce truncated variables below.
\[
X_{nk}^{(1)}  \buildrel \Delta \over =  - I(\omega _{nk} (x)\varepsilon _{nk}  \le  - 1) + \omega _{nk} (x)\varepsilon _{nk} I(\left| {\omega _{nk} (x)\varepsilon _{nk} } \right| < 1) + I(\omega _{nk} (x)\varepsilon _{nk}  \ge 1),
\]
\begin{center}
$X_{nk}^{}  \buildrel \Delta \over = \omega _{nk} (x)\varepsilon _{nk} I(\left| {\omega _{nk} (x)\varepsilon _{nk} } \right| < 1),$
and
$S_n (x) \buildrel \Delta \over = \sum\nolimits_{k = 1}^n {\omega _{nk} (x)\varepsilon _{nk} } ,S_n^{(1)} (x) \buildrel \Delta \over = \sum\nolimits_{k = 1}^n {X_{nk}^{(1)} } .$
\end{center}

From $E\left| X \right| < \infty$, it follows that
\[
nP\{ \left| X \right| > n\}  \to 0,_{}^{} as_{}^{}\  n \to \infty .
\]
Then for $\forall x \in C(g)$, when $n \to \infty$,
\[
P(S_n^{(1)} (x) \ne S_n (x)) \le \sum\limits_{k = 1}^n {P(\left| X \right| \ge \left| {\omega _{nk} (x)} \right|^{ - 1} )}  \to 0,
\]
since Condition $(B_2 )$ and $\mathop {\sup }\limits_k \left| {\omega _{nk} (x)} \right| = o(1)$.

It suffices to show that $S_n^{(1)} (x) $ converge to zero in probability for $\forall x \in C(g)$. Observe that $\{ X_{nk}^{(1)} ,1 \le k \le n\}$ are also NQD by lemma 2.1, hence by Chebyshev inequality,
\begin{eqnarray*}
\ \ \ \ \ \ \ P(\left| {S_n^{(1)} (x)} \right| > \varepsilon ) &\le& \varepsilon ^{ - 2} \sum\limits_{k = 1}^n {E(X_{nk}^{(1)} } )^2 I(\left| {\sum\limits_{k = 1}^n {X_{nk}^{(1)} } } \right| > \varepsilon )\\
&\le& \varepsilon ^{ - 2} [\sum\limits_{k = 1}^n {P(\left| {\omega _{nk} (x)\varepsilon _{nk} } \right| \ge 1)}  + \sum\limits_{k = 1}^n {\omega _{nk}^2 (x)E\varepsilon _{nk}^2 I(\left| {\omega _{nk} (x)\varepsilon _{nk} } \right| < 1)} ]\\
 &\buildrel \Delta \over =& \varepsilon ^{ - 2} (I_{n1}  + I_{n2} ), \ \ \ \ \ \ \ \ \ \ \ \ \ \ \ \ \ \ \ \ \ \ \ \ \ \ \ \ \ \ \ \ \ \ \ \ \ \ \ \ \ \ \ \ \ \ \ \ \ \ \ \ \ \ \ \ \ \  (4.5)
\end{eqnarray*}
where the first inequality is due to lemma 2.2.

When it come to $I_{n1}$,
\[
I_{n1}  \le \sum\limits_{k = 1}^n {P(\left| {\omega _{nk} (x)X} \right| \ge 1)}  \le \mathop {\sup }\limits_k \left| {\omega _{nk} (x)} \right|\sum\limits_{k = 1}^n {E\left| X \right|I(\left| {\omega _{nk} (x)X} \right| \ge 1)}  \to 0,_{}  \eqno (4.6)
\]
as $ n \to \infty $.

Now, choose a number $A$ such that for $x \ge A$, $P(\left| X \right| \ge x) \le \varepsilon x^{ - 1}$ .
We have
\begin{eqnarray*}
t^{ - 1} \int_0^t {xP(\left| {\varepsilon _k } \right| \ge x)} dx &\le& t^{ - 1} \int_0^t {xP(\left| X \right| \ge x)} dx \le t^{ - 1} M + \varepsilon ,
\end{eqnarray*}
which means that $t^{ - 1} \int_0^t {xP(\left| {\varepsilon _k } \right| \ge x)} dx \to 0$ when $t \to \infty$.

  Again
\[
0 \le \int_{\left| x \right| < t}^{} {x^2 } dP(\varepsilon _k  < x) =  - t^2 P(\left| {\varepsilon _k } \right| \ge t) + 2\int_0^t {xP(\left| {\varepsilon _k } \right| \ge x)} dx \le 2\int_0^t {xP(\left| X \right| \ge x)} dx.
\]

Therefore, it follows that
\begin{eqnarray*}
\ \ \ \ \ \ \ \ \ \ \ \ \ \ \ \ \ \ I_{n2}&=&\sum\limits_{k = 1}^n {\omega _{nk}^2 (x)E\varepsilon _k^2 I(\left| {\omega _{nk} (x)\varepsilon _k } \right| < 1)}\\
 &\le& 2\sum\limits_{k = 1}^n {\omega _{nk}^2 (x)} \int_0^{\left| {\omega _{nk} (x)} \right|^{ - 1} } {uP(\left| {\varepsilon _k } \right| \ge u)du}\\
 &\le& 2\sum\limits_{k = 1}^n {\omega _{nk}^2 (x)\int_0^{\left| {\omega _{nk} (x)} \right|^{ - 1} } {uP(X \ge u)du} }  \to 0,_{} \ as_{}^{}\  n \to \infty .   \ \ \ \ \ \ \ \ \ \ \ \ (4.7)
 \end{eqnarray*}

 Hence, the theory follows from $(4.4) \sim (4.7)$. This ends the proof.\ \ \ \ \ \ \ \ \ \ \ \ \ \ \ \ \ \ \qed \\
\textit{Proof of Corollary 3.2} \ By the discussion in Corollary 3.1, it is the direct result of Theorem 3.3. This completes the proof of Corollary 3.2.\ \ \ \ \ \ \ \ \ \ \ \ \ \ \ \ \ \ \ \ \ \ \ \ \ \ \ \ \ \ \ \ \ \ \ \ \ \ \ \ \ \ \ \ \ \ \ \  \ \ \ \ \ \ \ \ \ \ \ \ \ \ \ \ \ \qed \\
\\
\textbf{Competing interests}\\
The authors declare that they have no competing interests.

\small
\baselineskip=8pt

\end{document}